\documentclass{amsart}
\usepackage[all]{xy}

\usepackage{amsmath,amssymb,calligra,mathrsfs}
\usepackage[all]{xy}

\newtheorem{De}{Definition}[section]
\newtheorem{Th}[De]{Theorem}
\newtheorem{Prop}[De]{Proposition}
\newtheorem{Le}[De]{Lemma}
\newtheorem{Co}[De]{Corollary}

\def \Extalg{\mathop{\sf Extalg}\nolimits}

\numberwithin{equation}{subsection}

\begin{document}

\title{The cohomology of $\lambda$-rings and $K$-theory}

\author[Michael Robinson]{Michael Robinson}
\address{
Department of Mathematics\\
University of Leicester\\
University Road\\
Leicester\\
LE1 7RH, UK} \email{mr85@alumni.le.ac.uk}
%

\maketitle

\begin{abstract}
We introduce the Andr\'{e}-Quillen cohomology of $\lambda$-rings and $\Psi$-rings, this is different to the $\lambda$-ring cohomology defined by Yau in 2005. We show that there is a natural transformation connecting the cohomology of the $K$-theory of spheres to the homotopy groups of spheres.
\end{abstract}

\bigskip
\section{Introduction}
$\lambda$-rings were first introduced in an algebraic-geometry
setting by Grothendieck in 1958, then later used in group theory by
Atiyah and Tall. In 1962 Adams introduced the operations $\Psi^{i}$ to study vector fields of spheres. These operations give us another type of ring, the $\Psi$-rings, which are closely related to the $\lambda$-rings. The main example of $\lambda$-rings and $\Psi$-rings are in algebraic topology; the K-theory of a topological space is a $\lambda$-ring and $\Psi$-ring.

For more detailed information on the cohomology of $\lambda$-rings and $\Psi$-rings, see my thesis \cite{Thesis}. In 2005, Donald Yau defined a cohomology for $\lambda$-rings. We are using the Andr\'{e}-Quillen cohomology of $\lambda$-rings and $\Psi$-rings which is different to Yau's cohomology. In this paper, we let $\mathbb{N}$ denote the non-zero natural numbers and $\mathbb{N}_{0}$ denote $\mathbb{N}\cup\{0\}$.

\section{$\lambda$-rings and $\Psi$-rings}
In this section, we introduce the definitions of a $\lambda$-ring and a $\Psi$-ring. For more information on $\lambda$-rings and $\Psi$-rings, see Atiyah and Tall \cite{AT} or Knutson \cite{Knutson}.

\subsection{$\lambda$-rings}
A \textit{$\lambda$-ring} is a unital commutative ring, $R$, together with a sequence of operations $\lambda^{i}: R \rightarrow
R$, for $i \in \mathbb{N}_{0}$, satisfying
\begin{enumerate}
    \item $\lambda^{0}(r) = 1,$
    \item $\lambda^{1}(r) = r,$
    \item $\lambda^{i}(r+s)= \Sigma_{k=0}^{i} \lambda^{k}(r)\lambda^{i-k}(s),$
    \item $\lambda^{i}(r) = 0$ for $i>1$,
    \item $\lambda^{i}(rs) = P_{i}(\lambda^{1}(r),\lambda^{2}(r),\ldots,\lambda^{i}(r),\lambda^{1}(s),\ldots,\lambda^{i}(s))$,
    \item $\lambda^{i}(\lambda^{j}(r)) =
    P_{i,j}(\lambda^{1}(r),\ldots,\lambda^{ij}(r)),$
\end{enumerate}
where $P_{i}$ and $P_{i,j}$ are universal polynomials with integer coefficients, see the cited material for precise definitions. Note that what we refer to as a $\lambda$-ring is called a special $\lambda$-ring in the materials.

\subsection{$\lambda$-modules}
$M$ is a \textit{$\lambda$-module} over the $\lambda$-ring $R$ if
$M$ is an $R$-module together with a sequence of group homomorphisms
$\Lambda^{i}: M \rightarrow M$, for $i \in \mathbb{N}$, satisfying
\begin{enumerate}
\item $\Lambda^{1}(m) = m,$
\item $\Lambda^{i}(rm) = \Psi^{i}(r) \Lambda^{i}(m),$
\item $\Lambda^{ij}(m) = (-1)^{(i+1)(j+1)}\Lambda^{i}\Lambda^{j}(m),$
\end{enumerate}
for all $m \in M, r \in R$ and $i,j \in \mathbb{N}$.

\subsection{$\lambda$-derivations}
A \textit{$\lambda$-derivation} of $R$ with values in $M$ is an additive homomorphism $d:R \rightarrow M$ such that
\begin{enumerate}
\item $d(r s) = r d(s) + d(r)s,$
\item $d(\lambda^{i}(r)) = \Lambda^{i}(d(r)) + \sum_{j=1}^{i-1}\Lambda^{j}(d(r))\lambda^{i-j}(r),$
\end{enumerate}
for all $r,s \in R$, and $i \in \mathbb{N}$. We let
$Der_{\lambda}(R,M)$ denote the set of all $\lambda$-derivations of
$R$ with values in $M$.

\subsection{$\Psi$-rings}
The $\lambda$-operations are neither additive nor multiplicative which makes them difficult to use. From these we can obtain the Adams operations
which are ring homomorphisms. We define the Adams operations using the Newton formula.
\[
\Psi^{i}(r) -
\lambda^{1}(r)\Psi^{i-1}(r)+...+(-1)^{i-1}\lambda^{i-1}(r)\Psi^{1}(r)
+ (-1)^{i}i\lambda^{i}(r) = 0.
\]

A \textit{$\Psi$-ring} is a unital commutative ring, $R$, together
with a sequence of ring homomorphisms $\Psi^{i}: R \rightarrow R$,
for $i \in \mathbb{N}$, satisfying
\begin{enumerate}
    \item $\Psi^{1}(r) = r,$
    \item $\Psi^{i}( \Psi^{j}(r)) = \Psi^{ij}(r),$
\end{enumerate}
for all $r \in R$ and $i,j \in \mathbb{N}$. We say that a $\Psi$-ring $R$ is \textit{special} if it also satisfies the property
\[\Psi^{p}(r) \equiv r^{p} \textrm{\qquad mod $pR$} \]
for all primes $p$ and $r \in R$. All of the $\Psi$-rings which come
from $\lambda$-rings are special. Wilkerson \cite{wilkerson} gives
us a condition for when the converse is true.

\begin{Th}\label{wilkth}(Wilkerson)
If $R$ is a $\mathbb{Z}$ torsion free special $\Psi$-ring, then
there exists a unique $\lambda$-ring structure on $R$ whose adams
operations are precisely the $\Psi$-operations.
\end{Th}

There is a unique $\lambda$-ring structure on the ring of the integers $\mathbb{Z}$ given by \[\lambda^{i}(x) = \dbinom{x}{i},\] for $i \in \mathbb{N}$. The corresponding special $\Psi$-operations on $\mathbb{Z}$ are given by $\Psi^{i}(x) = x$ for $i \in \mathbb{N}$.

The definition of the free $\lambda$-ring is well-known, see
\cite{Knutson}. We are now going to construct the free $\Psi$-ring
on one generator $a$. Let $A$ be the free commutative ring generated
by $\{a_{i} | i \in \mathbb{N} \}$. Let the operations $\Psi^{i}: A
\rightarrow A$ be given by $\Psi^{i}(a_{j})=a_{ij}$, for $i,j \in
\mathbb{N}$. Then $A$ is the \textit{free $\Psi$-ring on one
generator}.

\begin{Le}
If $R$ and $S$ are $\Psi$-rings, then $R \otimes S$ with $\Psi^{i}:R \otimes S \rightarrow R \otimes S$ given by $\Psi^{i}(r,s) = (\Psi^{i}(r),\Psi^{i}(s)) $ is the coproduct in
the category of $\Psi$-rings.
\end{Le}

\begin{proof}
The coproduct of two commutative rings is given by the tensor
product, so we only need to check the $\Psi$-operations. There is a
unique $\Psi$-ring structure on $R \otimes S$ such that
\[ R \rightarrow R \otimes S, \qquad r \mapsto r \otimes 1,\]
\[ S \rightarrow R \otimes S, \qquad s \mapsto 1 \otimes s,\]
are homomorphisms of $\Psi$-rings given by
\begin{align*}\Psi^{i}(r \otimes s) &=\Psi^{i}((r \otimes 1)(1 \otimes s)) \\&= \Psi^{i}(r \otimes 1)\Psi^{i}(1 \otimes s) \\ &= (\Psi^{i}(r) \otimes 1)(1 \otimes \Psi^{i}(s)) \\&= \Psi^{i}(r) \otimes \Psi^{i}(s).
\end{align*}
\end{proof}

\begin{Co}
Let $A$ be the free commutative ring generated by
$\{a_{i},b_{i},\ldots,x_{i} | i \in \mathbb{N}\}$.
Let the operations $\Psi^{i}: A \rightarrow A$
be given by $\Psi^{i}(a_{j})=a_{ij}$, $\Psi^{i}(b_{j})=b_{ij}$, $\ldots$, $\Psi^{i}(x_{j})=x_{ij}$ for $i,j \in \mathbb{N}$. Then $A$ is
the \textit{free $\Psi$-ring} generated by $\{a, b, \ldots, x\}$.
\end{Co}

\subsection{$\Psi$-modules}
$M$ is a \textit{$\Psi$-module} over the $\Psi$-ring $R$ if $M$ is an
$R$-module together with a sequence of group homomorphisms $\psi^{i}:M
\rightarrow M$, for $i \in \mathbb{N}$, satisfying
\begin{enumerate}
    \item $\psi^{1}(m) = m,$
    \item $\psi^{i}(rm) = \Psi^{i}(r)\psi^{i}(m),$
    \item $\psi^{i}(\psi^{j}(m))  = \psi^{ij}(m),$
\end{enumerate}
for all $m \in M$, $r \in R$, and $i,j \in \mathbb{N}$. We let $R\mathfrak{-mod}_{\Psi}$ denote the category of all $\Psi$-modules over $R$.
We say that $M$ is \textit{special} if $R$ is special and
\[\psi^{p}(m) \equiv 0 \textrm{\qquad mod $pM$} \]
for all primes $p$ and $m \in M$.

\subsection{$\Psi$-derivations}
A \textit{$\Psi$-derivation} of $R$ with values in $M$ is an additive
homomorphism $d:R \rightarrow M$ such that
\begin{enumerate}
\item $d(r s) = r d(s) + d(r)s,$
\item $\psi^{i}(d(r)) = d(\Psi^{i}(r)),$
\end{enumerate}
for all $r, s \in R,$ and $i \in \mathbb{N}$. We let
$Der_{\Psi}(R,M)$ denote the set of all $\Psi$-derivations of $R$
with values in $M$.

\section{Cohomology of $\lambda$-rings}
It is known that there is an adjoint pair of functors
\[\xymatrix{\mathfrak{Sets}  \ar@<.5ex>[rr]^{F} && \lambda-\mathfrak{rings} \ar@<.5ex>[ll]^{U}}\]
where $U$ is the forgetful functor and $F$ takes a set $S$ to the
free $\lambda$-ring generated by $S$. The adjoint pair gives rise to
a comonad $\mathbb{G}$ on $\lambda-\mathfrak{rings}$ which is
monadic. Let $R$ be a $\lambda$-ring and $M$ be a $\lambda$-module
over $R$. We define the cohomology of the $\lambda$-ring $R$ with
coefficients in $M$, denoted by $H^{*}_{\lambda}(R,M)$, to be the comonad cohomology \cite{BB} of
$R$ with coefficients in $Der_{\lambda}(-,M)$. Note that $Der_{\lambda}(-,M)$ is a functor from the category of $\lambda$-rings to the category of abelian groups.

\begin{Co} For any $\lambda$-ring $R$ and $M \in R\mathfrak{-mod}_{\lambda}$, we have
\[H^{0}_{\lambda}(R,M) \cong Der_{\lambda}(R,M).\]
Furthermore, if $R$ is free as a $\lambda$-ring then $H^{i}_{\lambda}(R,M) =0$ for $i>0$.
\end{Co}

Let $R$ be a $\lambda$-ring and $M \in R\mathfrak{-mod}_{\lambda}$. A \textit{$\lambda$-ring extension} of $R$ by $M$ is an exact
sequence
\[\xymatrix{0 \ar[r] & M \ar^{\alpha}[r] & X \ar^{\beta}[r] & R \ar[r] & 0}\] where
$X$ is a $\lambda$-ring, $\beta$ is a map of $\lambda$-rings,
$\alpha$ is an additive homomorphism such that $\alpha \lambda^{i} = \Lambda^{i} \alpha$ for all $i \in \mathbb{N}$ and $\alpha(m)x =
\alpha(m \beta(x))$ for all $m \in M$ and $x \in X$. The map $\alpha$ identifies $M$ with an ideal of square-zero
in $X$.

Two $\lambda$-ring extensions $(X),(X')$ with $R,M$ fixed are said
to be \textit{equivalent} if there exists a map of $\lambda$-rings
$\phi: X \rightarrow X'$ such that the following diagram commutes.
\[\xymatrix{0 \ar[r] & M \ar@{=}[d] \ar[r] & X \ar[r] \ar[d]^{\phi} & R \ar[r] \ar@{=}[d] & 0 \\ 0 \ar[r] & M \ar[r] & X' \ar[r] & R \ar[r] & 0}\]
We denote the set of equivalence classes of $\lambda$-ring
extensions of $R$ by $M$ by $Extalg_{\lambda}(R,M)$.

\begin{Le} For any $\lambda$-ring $R$ and $M \in R\mathfrak{-mod}_{\lambda}$, we have
\[H^{1}_{\lambda}(R,M) \cong Extalg_{\lambda}(R,M).\]
\end{Le}
The proof can be found in my thesis \cite{Thesis}.

\section{Cohomology of $\Psi$-rings}
Let $I$ denote  the category with one
object associated to the multiplicative monoid of the
nonzero natural numbers. We can consider
$\Psi$-rings as diagrams of commutative rings; $\Psi$-rings are
functors from $I$ to the category of commutative rings.
 \[ R: I \rightarrow \mathfrak{Com.rings.}\]
It is well known that there is an adjoint pair of functors
\[\xymatrix{\mathfrak{Sets}  \ar@<.5ex>[rr]^{F} && \mathfrak{Com.rings} \ar@<.5ex>[ll]^{U}}\]
where $U$ is the forgetful functor and $F$ takes a set $S$ to the free commutative ring generated by $S$. The adjoint pair gives rise to a comonad $\mathbb{G}$ on $\mathfrak{Com.rings}$ which is monadic and the cohomology with respect to this comonad is the Andr\'{e}-Quillen cohomology of commutative rings. The adjoint pair gives rise to another adjoint pair
\[\xymatrix{\mathfrak{Sets}  \ar@<.5ex>[rr]^{F_{I}} && \mathfrak{Com.rings}^{I} \ar@<.5ex>[ll]^{U_{I}}}\]
where $U_{I}$ is the forgetful functor and $F_{I}$ takes a set $S$ to the free $\Psi$-ring generated by $S$.
This adjoint pair yields a comonad $\mathbb{G}_{I}$ on $\mathfrak{Com.rings}^{I} = \Psi-\mathfrak{rings}$ which is monadic.

Let $R$ be a $\Psi$-ring and $M$ be a $\Psi$-module over $R$. We define the cohomology of a $\Psi$-ring $R$ with coefficients in $M$, denoted by $H^{*}_{\Psi}(R,M)$, to be the comonad cohomology of $R$ with coefficients in $Der_{\Psi}(-,M)$. Note that $Der_{\Psi}(-,M)$ is a functor from the category of $\Psi$-rings to the category of abelian groups.

\begin{Co} For any $\Psi$-ring $R$ and $M \in R\mathfrak{-mod}_{\Psi}$, we have
\[H^{0}_{\Psi}(R,M) \cong Der_{\Psi}(R,M).\]
Furthermore, if $R$ is free as a $\Psi$-ring then $H^{i}_{\Psi}(R,M) =0$ for $i>0$.
\end{Co}

Let $R$ be a $\Psi$-ring and $M \in R\mathfrak{-mod}_{\Psi}$. A
\textit{$\Psi$-ring extension} of $R$ by $M$ is an exact sequence
\[\xymatrix{0 \ar[r] & M \ar^{\alpha}[r] & X \ar^{\beta}[r] & R \ar[r] & 0}\] where
$X$ is a $\Psi$-ring, $\beta$ is a map of $\Psi$-rings, $\alpha$ is an additive homomorphism such that $\alpha \Psi^{i} = \psi^{i} \alpha$ for all $i \in \mathbb{N}$ and $\alpha(m)x =
\alpha(m \beta(x))$ for all $m \in M$ and $x \in X$. The map $\alpha$ identifies $M$ with an ideal of square-zero
in $X$.

Two $\Psi$-ring extensions $(X),(\overline{X})$ with $R,M$ fixed are
said to be \textit{equivalent} if there exists a map of $\Psi$-rings
$\phi: X \rightarrow \overline{X}$ such that the following diagram
commutes.
\[\xymatrix{0 \ar[r] & M \ar@{=}[d] \ar[r] & X \ar[r] \ar[d]^{\phi} & R \ar[r] \ar@{=}[d] & 0 \\ 0 \ar[r] & M \ar[r] & \overline{X} \ar[r] & R \ar[r] & 0}\]
We denote the set of equivalence classes of $\Psi$-ring extensions
of $R$ by $M$ by $Extalg_{\Psi}(R,M)$. If $R$ and $M$ are special, then we say that an extension
\[\xymatrix{0 \ar[r] & M \ar^{\alpha}[r] & X \ar^{\beta}[r] & R \ar[r] & 0}\]
is \textit{special} if $X$ is also special.

\begin{Le} For any $\Psi$-ring $R$ and $M \in R\mathfrak{-mod}_{\Psi}$, we have
\[H^{1}_{\Psi}(R,M) \cong Extalg_{\Psi}(R,M).\]
\end{Le}

For each $n \in \mathbb{N}_{0}$, there is a natural system \cite{BW} on $I$ as follows
\[D_{f} := H_{AQ}^{n}(R,M^{f})\] where $M^{f}$ is an $R$-module with
$M$ as an abelian group with the following action of $R$
\[(r,m) \mapsto \Psi^{f}(r)m, \textrm{ for $r \in R, m \in M$}.\]
For any morphism $u \in I$, we have $u_{*}:D_{f} \rightarrow D_{uf}$ which is
induced by $\Psi^{u}: M^{f} \rightarrow M^{uf}$. For any morphism $v \in I$, we
have $v^{*}:D_{f} \rightarrow D_{fv}$ which is induced by $\Psi^{v}:
R \rightarrow R$.

\begin{Th} There exists a spectral sequence
\[E^{p,q}_{2} = H^{p}_{BW}(I,\mathcal{H}^{q}(R,M)) \Rightarrow H^{p+q}_{\Psi}(R,M)\]
where $\mathcal{H}^{q}(R,M)$ is the natural system on $I$ whose value on a morphism $\alpha$ in $I$
 is given by $H_{AQ}^{q}(R,M^{\alpha})$ and $H^{p}_{BW}(I,\mathcal{H}^{q}(R,M))$ is the Baues-Wirsching cohomology \cite{BW} of the small category $I$ with coefficients in the natural system $\mathcal{H}^{q}(R,M)$.
\end{Th}
The proof of this theorem can be found in my thesis \cite{Thesis} or the paper \cite{DiagCohom}.

\section{Natural transformation}
We let $K(-)$ denote the complex K-theory and $\widetilde{K}(-)$ denote the reduced complex K-theory. Let $X,Y$ be topological spaces such that $\widetilde{K}(Y)=0$  and $\widetilde{K}(\Sigma X) = 0$. Let $f:
Y \rightarrow X$ be a continuous map, then we can consider the Puppe sequence
$$\xymatrix{Y \ar^{f}[r] & X \ar[r] & C_{f} \ar[r] & \Sigma Y \ar[r] & \Sigma X \ar[r] & \Sigma C_{f} \ar[r] & \ldots}$$
where $C_{f}$ is the mapping cone of $f$, and $\Sigma X$ is the
suspension of $X$. After applying the functor $\widetilde{K}(-)$ we
get the long exact sequence.
$$\xymatrix{\ldots \ar[r] &  \widetilde{K}(\Sigma X) \ar[r] & \widetilde{K}(\Sigma Y) \ar[r] & \widetilde{K}(C_{f}) \ar[r] & \widetilde{K}(X) \ar[r] & \widetilde{K}(Y)}$$
However, since $\widetilde{K}(\Sigma X) = 0$ and
$\widetilde{K}(Y) = 0$ we obtain the short exact
sequence.
$$\xymatrix{ 0 \ar[r] & \widetilde{K}(\Sigma Y) \ar[r] & K(C_{f}) \ar[r] & K(X) \ar[r] & 0}$$
This gives us the following proposition.

\begin{Prop}
If $X$ and $Y$ are topological spaces as above then there exist
natural transformations \[\tau_{\lambda} : [Y,X] \rightarrow
Extalg_{\lambda}(K(X),\widetilde{K}(\Sigma Y)),\] \[\tau_{\Psi} : [Y,X]
\rightarrow Extalg_{\Psi}(K(X),\widetilde{K}(\Sigma Y)).\]
\end{Prop}

\begin{Co}
If $X$ is a topological space such that $\widetilde{K}(\Sigma X) = 0$
then there exist natural transformations $\tau_{\lambda} :
\pi_{2n-1}(X) \rightarrow Extalg_{\lambda}(K(X),\widetilde{K}(S^{2n}))$
and $\tau_{\Psi} : \pi_{2n-1}(X) \rightarrow
Extalg_{\Psi}(K(X),\widetilde{K}(S^{2n})).$
\end{Co}

\section{The Hopf invariant of an extension}
Consider the commutative ring $R$ generated by $x$ and $y$ as an abelian group, $R \cong \mathbb{Z}[x] \oplus \mathbb{Z}[y]$, where $x$ is the unit of the ring and $y^{2} = 0$. The ring $R$ is known as the ring of dual numbers. Let $M \cong \mathbb{Z}[z]$ be the $R$-module such that $y \cdot z = 0$. We can consider the extensions of $R$ by $M$ in the category of commutative rings.
All the extensions have the following form
\begin{equation}\label{aqext} \xymatrix{ 0 \ar[r] & M \ar[r]
& X \oplus \mathbb{Z}[\gamma] \ar[r] & R \ar[r] &0 }\end{equation} where $X \cong
\mathbb{Z}[\alpha] \oplus \mathbb{Z}[\beta]$ as an abelian group with $\alpha$ being the
image of the generator $z$, the image of the unit $\gamma$ is the unit $x$ and the image of
$\beta$ being the generator $y$. Since $M^{2}=0$ we get that $\alpha^{2}=0$. Since
$y^{2}=0$, we get that $\alpha \beta = 0$ and
$\beta^{2} = h \alpha$ for some integer $h$. We define $h$ to be the \textit{Hopf invariant} of the extension (\ref{aqext}).

We are going to consider the extensions of $K(S^{2n})$ by $\widetilde{K}(S^{2n'})$ in the category of $\Psi$-rings.
We are going to prove the following theorem
\begin{Th}\label{psiextttt}
 \[Extalg_{\Psi}(K(S^{2n}),\widetilde{K}(S^{2n'})) \cong  \left\{
                                                             \begin{array}{ll}
                                                               \mathbb{Z} \oplus
\mathbb{Z}_{G_{n,n'}} & \textrm{if } n\neq n'; \\
                                                               \mathbb{Z} \oplus \prod_{\textrm{p prime}} \mathbb{Z} & \textrm{if } n=n'.
                                                             \end{array}
                                                           \right.
   \]
where $G_{n,n'}$ denotes the greatest common divisor of all the integers in
the set \\$\{l^{n}-l^{n'} | l \in \mathbb{Z}, l\geq 2\}$
\end{Th}

\begin{Co}\label{spsiext}
If $n \neq n'$ then
\[Extalg_{\lambda}(K(S^{2n}),\widetilde{K}(S^{2n'})) \cong
\{(h,\nu) \in \mathbb{Z} \oplus \mathbb{Z}_{G_{n,n'}}|h \equiv \nu
\frac{(2^{n}-2^{n'})}{G_{n,n'}} \textrm{ mod }2.\}\] If $n = n'$
then
\begin{align*}Extalg_{\lambda}(K(S^{2n}),\widetilde{K}(S^{2n'}))
\cong \{(h,\nu_{2},\nu_{3},\ldots) \in \mathbb{Z} \oplus
\prod_{\textrm{p prime}} \mathbb{Z}|& h \equiv \nu_{2} \textrm{ mod
}2, \\ &\nu_{p}\equiv 0 \textrm{ mod }p, \textrm{
}p>2.\}\end{align*}
\end{Co}

All the $\Psi$-ring extensions of $K(S^{2n})$ by $\widetilde{K}(S^{2n'})$ have the form (\ref{aqext}).
The $\Psi$-operations on $\Psi^{k}:X \rightarrow X$ are given by
\[\psi^{k}(m,r)= (k^{n'}m + \nu_{k}r,k^{n}r)\]
for some $\nu_{k} \in \mathbb{Z}$.

$$\Psi^{k}(\Psi^{l}(m,r)) = (k^{n'}l^{n'}m + k^{n'}\nu_{l}r + \nu_{k}l^{n}r,k^{n}l^{n}r)$$
$$\Psi^{l}(\Psi^{k}(m,r)) = (l^{n'}k^{n'}m + l^{n'}\nu_{k}r + \nu_{l}k^{n}r,l^{n}k^{n}r)$$
Since the $\Psi$-operations commute, we get that
$$\nu_{l}r(k^{n'}-k^{n}) = \nu_{k}r(l^{n'}-l^{n})$$
If $n=n'$ then there is no restriction on the choice of $\nu_{p}$ for $p$ prime. Otherwise we can rearrange the above to get that
$$\nu_{l} = \nu_{k}\frac{(l^{n'}-l^{n})}{(k^{n'}-k^{n}).}$$
By setting $k=2$ we get that for all $l \geq 2$
$$\nu_{l} = \nu_{2}\frac{(l^{n'}-l^{n})}{(2^{n'}-2^{n}).}$$
We can write all the $\nu_{l}$'s as multiples of $\nu_{2}$ since
$$\nu_{l} = \nu_{2}\frac{(l^{n'}-l^{n})}{(2^{n'}-2^{n})} = \nu_{2}\frac{(k^{n'}-k^{n})}{(2^{n'}-2^{n})}\frac{(l^{n'}-l^{n})}{(k^{n'}-k^{n})} = \nu_{k}\frac{(l^{n'}-l^{n})}{(k^{n'}-k^{n}).}$$

Since $\nu_{2}$ is an integer, we get that $\nu_{2} =
\frac{z(2^{n'}-2^{n})}{G_{n,n'}}$ for some integer z.

If we replace the generator $\beta$ by $\beta + N\alpha$, note that
$(\beta + N\alpha)^{2} = h\alpha$, then we have to replace $\nu_{k}$
by $\nu_{k}+N(k^{n'}-k^{n})$. We get that
$$\nu_{k}+N(k^{n'}-k^{n}) =
\nu_{2}\frac{k^{n'}-k^{n}}{2^{n'}-2^{n}} +N(k^{n'}-k^{n}) =
\frac{(\nu_{2}+N(2^{n'}-2^{n}))(k^{n'}-k^{n})}{(2^{n'}-2^{n})}
$$
So we only have to be concerned with replacing $\nu_{2}$ by
$\nu_{2}+N(2^{n'}-2^{n})$, then our usual formula for $\nu_{k}$
holds. Hence we are replacing $\frac{z(2^{n'}-2^{n})}{G_{n,n'}}$ by
$$\frac{z(2^{n'}-2^{n})}{G_{n,n'}} + N(2^{n'}-2^{n}) = \frac{(z+NG_{n,n'})(2^{n'}-2^{n})}{G_{n,n'}}$$

This proves theorem \ref{psiextttt}. The isomorphism depends on $n$
and $n'$. By restricting to the special $\Psi$-ring extensions, we
get that $\nu_{2}r \equiv hr^{2}$ mod $2$ and $\nu_{p}r \equiv 0$
mod $p$ for $p \geq 3$. Since all the $\Psi$-rings in our extensions
are $\mathbb{Z}$ torsion free, the theorem of Wilkerson \ref{wilkth}
gives us corollary \ref{spsiext}.

\begin{Prop}
If there exists an extension in
$Extalg_{\lambda}(K(S^{2n}),\widetilde{K}(S^{2n'}))$ whose Hopf
invariant is odd, then either $n=n'$ or $min(n,n') \leq
g^{2}_{|n-n'|}$, where $g^{p}_{j}$ denotes the multiplicity of the
prime p in the prime factorisation of the greatest common divisor of
the set of integers $\{(k^{j}-1 ) |\textrm{ } k \in
\mathbb{N}-\{1,qp | \forall q \in \mathbb{N} \}\}.$
\end{Prop}

\begin{proof}The case when $n=n'$ is clear. Assume that $n \neq n'$, then the special $\Psi$-ring extensions are given by a pair $(h,\nu)$ where $h$ is the Hopf invariant. By \ref{spsiext}, $h$ can only be odd if $2^{n}$ divides $G_{n,n'}$. Assume that $n<n'$, since the other case is analogous. The multiplicity of 2 in the prime factorisation of $G_{n,n'}$ is $n$ if $n \leq g^{2}_{|n-n'|}$ or $g^{2}_{|n-n'|}$ if $g^{2}_{|n-n'|} < n$. It follows that if $n \leq g^{2}_{|n-n'|}$ then $2^{n}$ divides $G_{n,n'}$.
\end{proof}

Note that $g^{2}_{2n-1}=1$ for all $n\in \mathbb{N}$. Since $(k^{2n}-1)=(k^{n}+1)(k^{n}-1)$ it follows that $g^{2}_{2n}=\left\{
\begin{array}{ll}
3, & n \textrm{ odd} \\
g^{2}_{n}+1, & n \textrm{ even.}
\end{array}
\right.
$

\begin{Th}
If there exists an extension in
$Extalg_{\lambda}(K(S^{2n}),\widetilde{K}(S^{2n'}))$ whose Hopf
invariant is odd, then one of the following is satisfied
\begin{enumerate}
  \item $n=n'$.
  \item $n=1$ or $n'=1$.
  \item $n'-n$ is even and either $n=2$ or $n'=2$.
  \item $n'>n\geq 3$ and $n' = n+2^{n-2}b$ for some $b \in \mathbb{N}_{0}$.
  \item $n>n'\geq 3$ and $n = n'+2^{n'-2}b$ for some $b \in \mathbb{N}_{0}$.
\end{enumerate}
\end{Th}

\begin{proof}
1. is clear. \\2. follows from $g^{2}_{n} \geq 1$ for all
$\mathbb{N}$.\\3. follows from $g^{2}_{2n} \geq 3$ for all $n \in
\mathbb{N}$.
\\4. and 5. follows from $g^{2}_{|n-n'|}$ being 2 plus the multiplicity of 2 in the prime factorisation of $|n-n'|$.
\end{proof}

\begin{Le}
If there exists an extension in
$Extalg_{\lambda}(K(S^{2n}),\widetilde{K}(S^{2an}))$ for $a \in
\mathbb{N}$ whose Hopf invariant is odd, then one of the following
is satisfied
\begin{enumerate}
  \item $n=1,2$ or $4$.
  \item $n=3$ and $a$ is even.
  \item $n \geq 5$ and $(a-1)n = n + 2^{n-2}b$ for some $b \in \mathbb{N}_{0}$.
\end{enumerate}
\end{Le}

\begin{Co}
If there exists an extension in
$Extalg_{\lambda}(K(S^{2n}),\widetilde{K}(S^{4n}))$ whose Hopf
invariant is odd, then $n=1,2$ or $4$.
\end{Co}

\begin{Co}[Adams]
If $f: S^{4n-1} \rightarrow S^{2n}$ is a continuous map whose Hopf invariant is odd, then $n=1,2$ or $4$.
\end{Co}

\section{Stable Extalg groups of spheres}
\begin{Prop}
If $n> k+1$ then $G_{n,n+k} = G_{n+1,n+k+1}$.
\end{Prop}

\begin{proof}
Let $n>k+1$. We know that $G_{n,n+k} = G_{n+1,n+k+1}$ if and only if the multiplicity of any prime $p$ in the prime factorization of $G_{n,n+k}$ is $g^{p}_{k}$. For all primes $p>2$ we get that $p^{n} > 2^{k}-1$, so the multiplicity of $p$ in the prime factorisation of $G_{n,n+k}$ is $g^{p}_{k}$. We can easily see that $g^{2}_{k} \leq k+1$ for all $k$. It follows that the multiplicity of $2$ in the prime factorisation of $G_{n,n+k}$ is $g^{2}_{k}$.
\end{proof}

\begin{Co}
If $n> k+1$ then
\[\Extalg_{\lambda}(K(S^{2n}),\widetilde{K}(S^{2(n+k)})) \cong
\Extalg_{\lambda}(K(S^{2(n+1)}),\widetilde{K}(S^{2(n+k+1)})).\]
\end{Co}

The groups $\Extalg_{\lambda}(K(S^{2n}),\widetilde{K}(S^{2(n+k)}))$
are independent of $n$ for $n>k+1$, we call these the \textit{stable
Extalg groups of spheres} which we denote by $\Extalg^{s}_{2k}$.

\begin{Prop}
There are natural transformations
\[\Upsilon_{k}: \pi^{s}_{2k-1} \rightarrow \Extalg^{s}_{2k} \]
where $\pi^{s}_{2k-1}$ denotes the stable homotopy groups of spheres.
\end{Prop}
For small k these groups look as follows.
\begin{center}
\begin{tabular}{|l | l | l |}
  \hline
  k & $\pi^{s}_{2k-1}$ & $\Extalg^{s}_{2k}$ \\ \hline
  1 & $\mathbb{Z}_{2}$ & $2\mathbb{Z} \oplus \mathbb{Z}_{2}$ \\
  2 & $\mathbb{Z}_{24}\oplus \mathbb{Z}_{3}$ & $2\mathbb{Z}  \oplus \mathbb{Z}_{24}$ \\
  3 & $0$               & $2\mathbb{Z}\oplus \mathbb{Z}_{2}$ \\
  4 & $\mathbb{Z}_{240}$ & $2\mathbb{Z}\oplus \mathbb{Z}_{240}$ \\
  5 & $\mathbb{Z}_{2} \oplus \mathbb{Z}_{2} \oplus \mathbb{Z}_{2}$ & $2\mathbb{Z} \oplus \mathbb{Z}_{2}$ \\
  6 & $\mathbb{Z}_{504}$ & $2\mathbb{Z}\oplus \mathbb{Z}_{504}$ \\
  7 & $\mathbb{Z}_{3}$ & $2\mathbb{Z}\oplus \mathbb{Z}_{2}$ \\
  8 & $\mathbb{Z}_{480}\oplus \mathbb{Z}_{2}$ & $2\mathbb{Z} \oplus \mathbb{Z}_{480}$ \\
  \hline
\end{tabular}
 \end{center}


\begin{thebibliography}{999999999}

\bibitem{AT} {\sc M.F. Atiyah and D.O. Tall}. Group
representations, $\lambda$-rings and the $J$-homomorphism, Topology
8, 1969. p253-297.



\bibitem{BW} {\sc H.J. Baues and G. Wirsching}. Cohomology of small categories, Journal of pure and applied algebra 38, 1984.

\bibitem{BB} {\sc J.M. Beck}. Triples, algebras and cohomology, Ph.D. thesis, Columbia University,
1967.


\bibitem{Knutson} {\sc D. Knutson}. $\lambda$-Rings and the Representation Theory of the Symmetric Group, Springer, 1973, vol. 308.


\bibitem{Thesis} {\sc M. Robinson}.  The cohomology of $\lambda$-rings and $\Psi$-rings, Ph.D. thesis, University of Leicester, 2010.




\bibitem{DiagCohom} {\sc M. Robinson}. Cohomology of diagrams of algebras. arXiv:0802.3651v1 [math.KT].

\bibitem{Whitehead} {\sc G. Whitehead}.  Recent Advances in Homotopy Theory, The MIT Press, 1971.

\bibitem{wilkerson} {\sc C. Wilkerson}. Lambda rings, binomial
domains and vector bundles over $CP()$, Comm. Algebra 10 (1982),
311-328.


\bibitem{Yau} {\sc D. Yau}. Cohomology of $\lambda$-rings. J. Algebra 284 (2005), 37-51.

\end{thebibliography}
\end{document}